# 4x4 Matrix Representation of Hybrid Numbers, Gram Matrix and Golden Ratio

Abdullah MAĞDEN

Bursa Technical University

abdullah.magden@btu.edu.tr

**Abstract**

This study investigates the Lie algebra structure and geometric properties of hybrid numbers, under a specific non-associative multiplication table. Within the scope of the study, the structure constants of the system were derived through commutator relations, and 4x4 matrix representation was constructed. By defining the Gram matrix, complex and hybrid conjugates, the connection of the eigenvalues of Gram matrices to the golden ratio has been obtained. Theoretical analysis reveal that this hybrid Lie algebra possesses a semi-simple structure and that the characteristic spectrum of its Killing form is directly related to the Golden Ratio ($\varphi$). Furthermore, it is proven that the Complex and Hyperbolic conjugation operations are symmetry transformations that leave the Killing form invariant.

## 1. Introduction

Complex, hyperbolic and dual numbers are well known two dimensional number systems. Especially in the last century, a lot of researchers deal with the geometric and physical applications of these numbers [1,5,7].

Dual numbers were introduced in 1873 by William Clifford. Dual numbers are extensively used in the quantum mechanics and classic mechanics of screws [3]. In this respect, our work shows similarities to quaternions, which hold a significant place in mathematics and physics. Quaternions have three imaginary units $(h^1, h^2, h^3)$ whose squares are negative, and their multiplication rules are non-commutative $( h^i \cdot h^j \neq h^j \cdot h^i)$ This property has made quaternions an ideal tool for modeling three-dimensional rotations, especially in computer graphics, robotics, and navigation.

The defined hybrid algebra also shares the property of being non-commutative, similar to quaternions. However, our structure offers a richer framework by allowing the squares of its imaginary units to take different values $h^1 \cdot h^1 = -1$, $h^2 \cdot h^2 = 0$, $h^3 \cdot h^3 = 1$ [1]. This indicates that our algebra has the potential to model various transformations and physical phenomena, particularly wave propagation, in addition to rotations. This article provides a new perspective in abstract algebra and investigates how this structure can be applied to physical equations.

This article explores a novel algebraic structure for hybrid numbers $Z = a + b_i h^i$ that unifies the standard systems of complex, dual, and hyperbolic numbers. Three imaginary units $(h^1, h^2, h^3,)$ this algebra deviates from traditional systems due to its non-commutative multiplication. This structure not only offers a new perspective in abstract algebra but also holds the potential to model physical phenomena, particularly wave propagation.

## 2. Definition and Properties of the Hybrid Algebra

The set of hybrid numbers will be denoted by $\mathbb{H}$ and contains complex, dual and hyperbolic numbers as well as combined and mixed states of these types of three numbers. The corresponding geometry of hybrid numbers will be most general geometry including their combination as well as the Euclidean, Minkowski, and Galilean plane geometries.

Now, we define hybrit numbers and define fundemental multiplication rules for these units.

**Definition 2.1.** The set of hybrid numbers, denoted by $\mathbb{H}$ , is defined as [1]:

$$\mathbb{H} = \{a + b_i h^i : a, b_i \epsilon \mathbb{R}, i = 1,2,3\}$$

Here, $h^1 \cdot h^1 = -1$, $h^2 \cdot h^2 = 0$, $h^3 \cdot h^3 = 1$ and

| $\cdot$ | $1$ | $h^1$ | $h^2$ | $h^3$ | |
|---|---|---|---|---|---|
| $1$ | $1$ | $h^1$ | $h^2$ | $h^3$ | (2.1) |
| $h^1$ | $h^1$ | $-1$ | $1 - h^3$ | $h^2 + h^1$ | |
| $h^2$ | $h^2$ | $h^3 + 1$ | $0$ | $-h^2$ | |
| $h^3$ | $h^3$ | $-h^2 - h^1$ | $h^2$ | $1$ | |

Table 2.1

Here, we constructed using Einstein's summation convention and three imaginary units $(h^1, h^2, h^3)$ , this algebra deviates from traditional systems due to its non-commutative multiplication. The hybrid number system is a novel algebraic structure $Z = (a + b_i h^i)$ that unifies complex, dual, and hyperbolic numbers into a single framework. A defining characteristic of this algebra is its noncommutative multiplication, which distinguishes it from traditional 2D number systems.

$$h^1 \cdot h^3 = -h^3 \cdot h^1 = h^2 + h^1$$

Two hybrid numbers are equal if all their components are equal, one by one. The sum of two hybrid numbers is defined by summing their components. Addition operation in the hybrid numbers is both commutative and associative. Zero is the null element. With respect to the addition operation, the inverse element of $Z$ is $-Z$ , which is defined as having all the components of $-Z$ changed in their signs. This implies that, $(\mathbb{H}, +)$ is an Abelian group [1].

This property makes the algebra similar in structure to quaternions. The identity element of the algebra is defined as 1, but due to the nilpotent property of $h^2$ (its square is zero), there exist nonzero elements in this algebra that lack an inverse.

The product of two hybrid numbers $Z_1 \cdot Z_2$ involves a cross product of their vector parts $\overrightarrow{V_1} \times \overrightarrow{V_2}$. Because the vector cross product is noncommutative, $Z_1 \cdot Z_2 \neq Z_2 \cdot Z_1$. The fundamental multiplication rules for the three imaginary units $(h^1, h^2, h^3,)$ dictate that the order of multiplication determines the sign of the result.

**Definition 2.2. (Conjugate)**

$\bar{Z} = a - b_i h^i \epsilon \mathbb{H}$ is called the conjugate of $Z \epsilon \mathbb{H}$.

From (2.1) we write

$$\mathcal{C}(Z) = Z \cdot \bar{Z} = \left(a + b_i h^i\right)\left(a - b_i h^i\right) = (a+V)(a-V) = a^2 - V^2 \quad (2.2)$$
$$= a^2 + b_1^2 - 2b_1 b_2 - b_3^2, \qquad V = a + b_i h^i \ .$$

**Definition 2.3. (Character of Hybrid Number)**

$\mathcal{C}(Z) = Z \cdot \bar{Z} = a^2 + b_1^2 - 2b_1 b_2 - b_3^2 = a^2 + (b_1 - b_2)^2 - b_2^2 - b_3^2$ is called character of hybrid number $Z$[1].

We seen that from table 2.1. $Z \cdot \bar{Z} = \bar{Z} \cdot Z = a^2 + (b_1 - b_2)^2 - b_2^2 - b_3^2$.

This definition shows that hybrid numbers have a pseudo-norm structure.

In complex numbers, only the number 0 has a norm of zero. There are hybrid numbers whose $\mathcal{C}(Z)$ is zero even though they are not zero themselves. These numbers is called null hybrid numbers. If $\mathcal{C}(Z) = Z \cdot \bar{Z} = \bar{Z} \cdot Z = a^2 + b_1^2 - 2b_1 b_2 - b_3^2 = a^2 + (b_1 - b_2)^2 - b_2^2 - b_3^2 = 0$, than $Z$ is null hybrid number [1].

We can summarize the component structure of the norm with the following table:

| The component structure of the norm | | | |
|---|---|---|---|
| $+a^2$, real part $\geq 0$ | $+(b_1 - b_2)^2$, The difference $h^1 and\ h^2$ is always positive. | $-b_2^2$, dual part is always negative | $-b_3^2$ , hybrid part is always negative |

We can summarize the norm in subsystems as follows:

| Subsystem | Condition | $\mathcal{C}(Z)$ |
|---|---|---|
| Complex | $b_2 = b_3 = 0$ | $a^2 + b_1^2 \geq 0$ |
| Dual | $b_1 = b_3 = 0$ | $a^2 - b_2^2$ (uncertain) |
| Hyperbolic | $b_1 = b_2 = 0$ | $a^2 - b_3^2$ (Minkowskii) |
| $b_1 = b_2$ special situation | $b_1 - b_2 = 0$ | $-b_2^2 - b_3^2 \leq 0$ |
| Real | $b_1 = b_2 = b_3 = 0$ | $a^2 \geq 0$ |

## 3. 4x4 Matrix Representation of Hybrid Numbers

Amatrix representation for a hybrid number is especially important in order to facilitate multiplication of hybrid numbers. By defining an isomorphism between 4×4 matrices and hybrid numbers, we can easily multiply the hybrid numbers and prove many of their features. On the other hand, hybrid numbers can also be defined by considering the matrix representation. We can classify the hybrid numbers, with respect to determinant and discriminant of the characteristic equation of the 4 ×4 corresponding matrix.

Let be $Z = a + b_i h^i$ is a hybrid numbers. Multiplying the number Z by $h^1, h^2$ and $h^3$ respectively gives us the Left-Multiplication Matrix representation of the hybrid number as follows [5]:

$$M_L(Z) = \begin{bmatrix} a & -b_1 + b_2 & b_1 & b_3 \\ b_1 & a - b_3 & 0 & b_1 \\ b_2 & -b_3 & a + b_3 & b_1 - b_2 \\ b_3 & b_2 & -b_1 & a \end{bmatrix} \quad (3.1)$$

The left multiplication operator was used in the calculations. Similarly, Right-Multiplication Matrix the matrix can be found using the left product as follows [5]

$$M_R(Z) = \begin{bmatrix} a & -b_1 + b_2 & b_1 & b_3 \\ b_1 & a + b_3 & 0 & -b_1 \\ b_2 & b_3 & a - b_3 & -b_1 + b_2 \\ b_3 & -b_2 & b_1 & a \end{bmatrix} \quad (3.2)$$

If we use (3.1) find

$$M_L(I) = \begin{bmatrix} 1 & 0 & 0 & 0 \\ 0 & 1 & 0 & 0 \\ 0 & 0 & 1 & 0 \\ 0 & 0 & 0 & 1 \end{bmatrix} = I_4 , M_L(h^1) = \begin{bmatrix} 0 & -1 & 1 & 0 \\ 1 & 0 & 0 & 1 \\ 0 & 0 & 0 & 1 \\ 0 & 0 & -1 & 0 \end{bmatrix},$$

$$M_L(h^2) = \begin{bmatrix} 0 & 1 & 0 & 0 \\ 0 & 0 & 0 & 0 \\ 1 & 0 & 0 & -1 \\ 0 & 1 & 0 & 0 \end{bmatrix}, \ M_L(h^3) = \begin{bmatrix} 0 & 0 & 0 & 1 \\ 0 & -1 & 0 & 0 \\ 0 & -1 & 1 & 0 \\ 1 & 1 & 0 & 0 \end{bmatrix}.$$

From (2.1) we find that $[M_L(h^1)]^2 = -I_4$, $[M_L(h^2)]^2 = O_4$, $[M_L(h^3)]^2 = I_4$. If the operations in (2.1) are transferred to matrices, it is seen that all operations are validated with matrices. For example,

$$\begin{bmatrix} 0 & -1 & 1 & 0 \\ 1 & 0 & 0 & 1 \\ 0 & 0 & 0 & 1 \\ 0 & 0 & -1 & 0 \end{bmatrix} \cdot \begin{bmatrix} 0 & 1 & 0 & 0 \\ 0 & 0 & 0 & 0 \\ 1 & 0 & 0 & -1 \\ 0 & 1 & 0 & 0 \end{bmatrix} = \begin{bmatrix} 1 & 0 & 0 & -1 \\ 0 & 2 & 0 & 0 \\ 0 & 1 & 0 & 0 \\ -1 & 0 & 0 & 1 \end{bmatrix}$$

$$M_L(h^1) \cdot M_L(h^2) = M_L(1 - h^3)$$

$$\begin{bmatrix} 0 & 1 & 0 & 0 \\ 0 & 0 & 0 & 0 \\ 1 & 0 & 0 & -1 \\ 0 & 1 & 0 & 0 \end{bmatrix} \cdot \begin{bmatrix} 0 & -1 & 1 & 0 \\ 1 & 0 & 0 & 1 \\ 0 & 0 & 0 & 1 \\ 0 & 0 & -1 & 0 \end{bmatrix} = \begin{bmatrix} 1 & 0 & 0 & 1 \\ 0 & 0 & 0 & 0 \\ 0 & -1 & 2 & 0 \\ 1 & 0 & 0 & 1 \end{bmatrix}$$

$$M_L(h^2) \cdot M_L(h^1) = M_L(1 + h^3).$$

The others are obtained in a similar manner. From now on, we will use $M(Z)$ instead of $M_L(Z)$. Then we write a hybrid number $Z = a + b_i h^i$

$$M(Z) = aM(I) + b_iM(h^i) = aM(I) + b_1M(h^1) + b_2M(h^2) + b_3M(h^3)\,. \qquad (3.3)$$

**Theorem 3.1.** The map $\varphi: \mathbb{H} \to M_{4\times4}(\mathbb{R}), \varphi(Z) = M(Z)$ is a ring homomorphism.

**Proof.** Let be define the map $\varphi: \mathbb{H} \to M_{4\times4}(\mathbb{R}), \varphi(Z) = M(Z)$, where

$$M(Z) = M(a + b_ih^i) = \begin{bmatrix} a & -b_1 + b_2 & b_1 & b_3 \\ b_1 & a - b_3 & 0 & b_1 \\ b_2 & -b_3 & a + b_3 & b_1 - b_2 \\ b_3 & b_2 & -b_1 & a \end{bmatrix}, Z \in \mathbb{H}.$$

We will Show that this map is a ring homomorphism. Taking into account the addition and multiplication operation table defined for the hybrid numbers before, it can be easily shows.

**i.** $M(Z_1 + Z_2) = M(Z_1) + M(Z_2)$ ,

Let be $Z_1 = a + b_ih^i$, $Z_2 = c + d_ih^i \in \mathbb{H}$. We write $Z_1 + Z_2 = (a + c) + (b_i + d_i)h^i$. From definition of $M(Z)$, is obtained $M(Z_1 + Z_2) = M(Z_1) + M(Z_2)$.

**ii.** $M(Z_1 \cdot Z_2) = M(Z_1) \cdot M(Z_2)$, If we consider

$$M(Z_1) = \begin{bmatrix} a & -b_1 + b_2 & b_1 & b_3 \\ b_1 & a - b_3 & 0 & b_1 \\ b_2 & -b_3 & a + b_3 & b_1 - b_2 \\ b_3 & b_2 & -b_1 & a \end{bmatrix},$$

$$M(Z_2) = \begin{bmatrix} c & -d_1 + d_2 & d_1 & d_3 \\ d_1 & c - d_3 & 0 & d_1 \\ d_2 & -d_3 & c + d_3 & d_1 - d_2 \\ d_3 & d_2 & -b_1 & c \end{bmatrix}$$

then $M(Z_1 \cdot Z_2) = M(Z_1) \cdot M(Z_2)$ is obtained.

**iii.** $\varphi$ is injective. $\varphi(Z_1) = \varphi(Z_2) \Rightarrow M(Z_1) = M(Z_2) \Rightarrow Z_1 = Z_2$. Or , $M(Z) = O \Rightarrow a = 0, b_1 = 0, b_2 = 0, b_3 = 0$ $(from\ first\ column)$, that is, $Ker(\varphi)$ is trivial .

**Result 3.1.** The set $im(\varphi) \subset M_{4\times4}(\mathbb{R})$ forms an isomorphic subring with the hybrid number ring $\mathbb{H}$.

For it to be called an isomorphism, $\varphi$ would have to be surjective. However, $M_{4\times4}(\mathbb{R})$ is a much larger set. The correct term should be monomorphism or embedding. $\mathbb{H}$ is isomorphically embedded around 4x4 real matrices.

**Examble 3.1.** Let $Z_1 = 1 + h^1, Z_2 = h^2 + h^3$, $Z_1 \cdot Z_2 = 1 + h^1 + 2h^2$,

$$M(Z_1) \cdot M(Z_2) = \begin{bmatrix} a & -b_1 + b_2 & b_1 & b_3 \\ b_1 & a - b_3 & 0 & b_1 \\ b_2 & -b_3 & a + b_3 & b_1 - b_2 \\ b_3 & b_2 & -b_1 & a \end{bmatrix} \cdot \begin{bmatrix} c & -d_1 + d_2 & d_1 & d_3 \\ d_1 & c - d_3 & 0 & d_1 \\ d_2 & -d_3 & c + d_3 & d_1 - d_2 \\ d_3 & d_2 & -b_1 & c \end{bmatrix}$$

$$= \begin{bmatrix} 1 & -1 & 1 & 0 \\ 1 & 1 & 0 & 1 \\ 0 & 0 & 1 & 1 \\ 0 & 0 & -1 & 1 \end{bmatrix} \cdot \begin{bmatrix} 0 & 1 & 0 & 1 \\ 0 & -1 & 0 & 0 \\ 1 & -1 & 1 & -1 \\ 1 & 1 & 0 & 0 \end{bmatrix} = \begin{bmatrix} 1 & 1 & 1 & 0 \\ 1 & 1 & 0 & 1 \\ 2 & 0 & 1 & -1 \\ 0 & 2 & -1 & -1 \end{bmatrix} = M(Z_1 \cdot Z_2)$$

**3.2.Determinant of $M(Z)$**

$$Det\ M(Z) = \begin{vmatrix} a & -b_1 + b_2 & b_1 & b_3 \\ b_1 & a - b_3 & 0 & b_1 \\ b_2 & -b_3 & a + b_3 & b_1 - b_2 \\ b_3 & b_2 & -b_1 & a \end{vmatrix} = [a^2 + b_1^2 - 2b_1b_2 - b_3^2]^2 \quad (3.4)$$

$$= [a^2 + (b_1 - b_2)^2 - b_2^2 - b_3^2]^2$$

## Corollary 3.2.

- **i.** $Det\ M(Z) = [N(Z)]^2, DetM_R(Z) = DetM_L(Z)$ .
- **ii.** $trM(Z) = 4a$
- **iii.** $Det\ M(Z) \geq 0$, So, $M(Z)$ has an inverse matrix $M^{-1}(Z)$.

Since hybrid numbers are a combination of complex, dual, and hyperbolic numbers, only adverbs focusing on a specific subspace can be defined.

**Definition 2.4. (Complex conjugate, Hyperbolic conjugate and Anti conjugate)**

$\bar{Z}^{\mathbb{C}} = a - b_1h^1 + b_2h^2 + b_3h^3 \epsilon \mathbb{H}$ are called the complex conjugate,

$\bar{Z}^{\mathbb{H}} = a + b_1h^1 + b_2h^2 - b_3h^3 \epsilon \mathbb{H}$ the hyperbolic conjugate

and $\bar{Z}^{A} = -a + b_1h^1 + b_2h^2 + b_3h^3 \epsilon \mathbb{H}$ anti conjugate or scalar-vektor conjugate of $Z$.

Why is it appropriate to call it hyperbolic norm? This procedure is designed to analyze the hyperbolic (spacelike) orientation of the system. The $Z \cdot \bar{Z}^{\mathbb{H}}$-multiplier system highlights the Minkowski character. Real Part gives the Minkowski metric of the system. The $h^1, h^2$ components, however, don't eliminate the circular and dual components in the system through hyperbolic conjugation; rather, they interact with the hyperbolic component, moving it to a new plane. $h^3$ has completely disappeared, which indicates that the operation fully lives up to its name, hyperbolic conjugate.

Complex conjugate is used to analyze the effects of rotation in the Euclidean plane. Hyperbolic conjugate is used to study Lorentz transformations and hyperbolic geometries in the Minkowski plane. In the algebra of hybrid numbers, multiplying a number $Z$ by its various conjugate types is used to determine the geometric characteristic of the number and the subspace Euclidean, Minkowski, or Galilean to which it relates.

Let's decompose the number $Z$ into two parts: $A = a + b_1h^1 + b_2h^2$ circular-dual parts and a hyperbolic part $B = b_3h^3$. Than, we write $Z = A + B$ and $\bar{Z}^{\mathbb{H}} = A - B$. Therefore,

$Z \cdot \bar{Z}^{\mathbb{H}} = (A + B)(A - B) = A^2 - AB + BA - B^2$,

$$A^2 = (a + b_1h^1 + b_2h^2)^2 = a^2 + b_1^2(h^1)^2 + 2ab_1h^1 + 2ab_2h^2 + b_1b_2(h^1h^2 + h^2h^1)$$
$$= a^2 - b_1^2 + 2b_1b_2 + 2ab_1h^1 + 2ab_2h^2,$$
$$B^2 = (b_3h^3)^2 = b_3^2,$$
$$-AB = (a + b_1h^1 + b_2h^2)b_3h^3, \quad BA = (b_3h^3)(a + b_1h^1 + b_2h^2),$$
$$-AB + BA = -2b_1b_3h^1 + (2b_2b_3 - 2b_1b_3)h^2,$$

Then, we find

$$Z \cdot \bar{Z}^{\mathbb{H}} = \underbrace{(a^2 - b_1^2 + 2b_1b_2 - b_3^2)}_{Real\ part} + \underbrace{2(ab_1 - b_1b_3)h^1}_{h^1\ component} + \underbrace{2(ab_2 + b_2b_3 - b_1b_3)h^2}_{h^2\ component} \quad (3.5)$$

**Result 2.1.** $Z + \bar{Z}^A = 2V$.

Multiplications with other conjugates are found similarly. The coefficients of $h^1$, $h^2$ and $h^3$ in multiplication with all conjugates are obtained as follows:

| Conjugate | Real Part | $h^1$ | $h^2$ | $h^3$ |
|---|---|---|---|---|
| Normal | $a^2 + b_1^2 - 2b_1b_2 - b_3^2$ | 0 | 0 | 0 |
| Complex | $a^2 + b_1^2 + b_3^2$ | $2b_1b_3$ | $2ab_2 + 2b_1b_3$ | $2ab_3 - 2b_1b_2$ |
| Hyperbolic | $a^2 - b_1^2 + 2b_1b_2 - b_3^2$ | $2(ab_1 - b_1b_3)$ | $2(ab_2 + b_2b_3 - b_1b_3)$ | 0 |
| Anti | $-a^2 - b_1^2 + 2b_1b_2 + b_3^2$ | 0 | 0 | 0 |

Table 3.1

We called this matrix representation the standard left representation. The matrix representation for the conjugate is obtained as follows.

$$M(\bar{Z}^{\mathbb{H}}) = \begin{bmatrix} a & -b_1 + b_2 & b_1 & -b_3 \\ b_1 & a + b_3 & 0 & b_1 \\ b_2 & b_3 & a - b_3 & b_1 - b_2 \\ -b_3 & b_2 & -b_1 & a \end{bmatrix} \quad (3.6)$$

Matrix representations for other conjugates can be easily found.

**Remark 3.1** $Det\ M(Z) = DetM(\bar{Z}^{\mathbb{H}}),\ tr\ M(Z) = trM(\bar{Z}^{\mathbb{H}})$.

Multiple matrix representations can be written for hybrid numbers. These matrices are connected to each other via a similarity transformation. If we transform the generated matrix $M$ into the form $P \cdot M(Z) \cdot P^{-1}$ with an invertible matrix $P$, we obtain a new matrix representation. Although these matrices look different, the determinant and trace remain the same. The matrix representation of the hybrid number is obtained by expressing the algebraic multiplication operation as a linear operator. Under the transformation $L_Z: \mathbb{H} \to \mathbb{H}$, $L_Z(W) = Z \cdot W$, the coefficient matrix formed with respect to bases $\{1, h^1, h^2, h^3\}$ takes the form $M(Z)$. Different matrices can be created by choosing other bases. For example, if $\{1, h^1 + h^2, h^1 - h^2, h^3\}$were chosen as the base, other matrix representations could also be created.

**Theorem 3.1.** *All matrix representations defined for hybrid numbers, although varying depending on the chosen basis, preserve the algebraic structure. Any two valid matrix* *representations are connected by a similarity transformation, and under this transformation, the fundamental invariants of the algebra are preserved.[9]*

**Theorem 3.2.** *All matrix representations defined for a number $Z\epsilon\mathbb{H}$ are connected to each other by a matrix $P$. Under this transformation, Trace, Determinant, and Characteristic are polynomial invariants.[9]*

Multiple matrix representations can be written for hybrid numbers. These matrices are connected to each other via a similarity transformation. If we transform the generated matrix $M$ into the form $P \cdot M(Z) \cdot P^{-1}$ with an invertible matrix $P$, we obtain a new matrix representation. Although these matrices look different, the determinant and trace remain the same. The matrix representation of the hybrid number is obtained by expressing the algebraic multiplication operation as a linear operator[5,9]. Under the transformation $L_Z$: $\mathbb{H} \to \mathbb{H}$, $L_Z(W) = Z \cdot W$, the coefficient matrix formed with respect to bases $\{1, h^1, h^2, h^3\}$ takes the form $M(Z)$.

**Definition 3.1. (Gram Matrix)** *The matrix $G = (G_{ij})$ defined by*

$$G_{ij} = g(e_i, e_j) = Re(e_i \cdot \bar{e}_j) \tag{3.7}$$

*on some base is called the Gram matrix. Here, $\bar{e}_j$ represents the basis element to which the relevant conjugate operator has been applied, and*

$$Re(Z_1 \cdot Z_2) = a_1 a_1 + b_1 d_1 - b_1 d_2 - d_1 b_2 - b_3 d_3.$$

We find the Gram matrices for the conjugates defined above, according to the base $\{1, h^1, h^2, h^3\}$.

For $\bar{Z} = a - b_i h^i$, $\bar{e}_1 = 1, \bar{e}_2 = -h^1$, $\bar{e}_3 = -h^2$, $\bar{e}_3 = -h^3$ we find Gram matrice as,

$$G_O = \begin{bmatrix} 1 & 0 & 0 & 0 \\ 0 & 1 & -1 & 0 \\ 0 & -1 & 0 & 0 \\ 0 & 0 & 0 & -1 \end{bmatrix}, \tag{3.8}$$

Similarly for $\bar{Z}^{\mathbb{C}} = a - b_1 h^1 + b_2 h^2 + b_3 h^3 \epsilon\mathbb{H}$, $\bar{e}_1 = 1, \bar{e}_2 = -h^1$, $\bar{e}_3 = h^2$, $\bar{e}_3 = h^3$,

$$G_{\mathbb{C}} = \begin{bmatrix} 1 & 0 & 0 & 0 \\ 0 & 1 & 1 & 0 \\ 0 & -1 & 0 & 0 \\ 0 & 0 & 0 & 1 \end{bmatrix}, \tag{3.9}$$

$\bar{Z}^{\mathbb{H}} = a + b_1 h^1 + b_2 h^2 - b_3 h^3 \epsilon\mathbb{H}$, $\bar{e}_1 = 1, \bar{e}_2 = h^1$, $\bar{e}_3 = h^2$, $\bar{e}_3 = -h^3$,

$$G_{\mathbb{H}} = \begin{bmatrix} 1 & 0 & 0 & 0 \\ 0 & -1 & 1 & 0 \\ 0 & 1 & 0 & 0 \\ 0 & 0 & 0 & -1 \end{bmatrix} \tag{3.10}$$

$\bar{Z}^{A} = -a + b_1 h^1 + b_2 h^2 + b_3 h^3 \epsilon\mathbb{H}, \bar{e}_1 = -1, \bar{e}_2 = h^1$, $\bar{e}_3 = h^2$, $\bar{e}_3 = h^3$,

$$G_A = \begin{bmatrix} -1 & 0 & 0 & 0 \\ 0 & -1 & 1 & 0 \\ 0 & 1 & 0 & 0 \\ 0 & 0 & 0 & 1 \end{bmatrix}. \qquad (3.11)$$

**Corollary 3.1.** *For each non-degenerate symmetric Gram matrix,* $Det(G) = \mp 1$.

The Gram matrix $G_O,\ G_{\mathbb{H}}\ and\ G_A$ on the basis $\{1, h^1, h^2, h^3\}$ contains one of the following two blocks in terms of the $h^1 - h^2$ subblock as:

$$A = \begin{bmatrix} 1 & -1 \\ -1 & 0 \end{bmatrix},\ B = \begin{bmatrix} -1 & 1 \\ 1 & 0 \end{bmatrix}. \qquad (3.12)$$

Eigenvalues of both blocks are $\lambda^2 - \lambda - 1 = 0$ and $\lambda^2 + \lambda - 1 = 0$ respectively. Then we find eigenvalue of $A$

$$\lambda_1 = \varphi = \frac{1+\sqrt{5}}{2},\ \ (Golden\ Ratio),\ \lambda_2 = -\frac{1}{\varphi} = \frac{1-\sqrt{5}}{2}, \qquad (3.13)$$

and for $B$

$$\lambda_1 = -\varphi = -\frac{1+\sqrt{5}}{2},\ \ \lambda_2 = \frac{1}{\varphi} = \frac{-1+\sqrt{5}}{2}. \qquad (3.14)$$

**Corollary 3.2.**

i. *The eigenvalues of Gram matrices give the golden ratio.*
ii. *The eigenvalues for matrices A and B satisfy the following invariant equalities:*

$$|\lambda_1 \cdot \lambda_2| = 1, \qquad \lambda_1 + \lambda_2 = \mp 1, \qquad \frac{\lambda_1}{\lambda_2} = -\varphi^2.$$

The characteristic polynomial of block $h^1 - h^2$ satisfies equation

$$P(\lambda) = \lambda^2 \pm \lambda - 1 = 0$$

which is the defining equation of the golden ratio. The eigenvalues of the Gram matrices are given in the table below:

| Conjugate | Blok $h^1 - h^2$ | Eigenvalue of Gram Matrix | Det G | Signature |
|---|---|---|---|---|
| $\bar{Z}\ (conjugate)$ | $B$ | $\varphi, -\frac{1}{\varphi}, +1, -1$ | 1 | (2,2) |
| $\bar{Z}^{\mathbb{H}}\ (hyperbolik)$ | $A$ | $-\varphi, \frac{1}{\varphi}, +1, -1$ | 1 | (2,2) |
| $\bar{Z}^{A}\ (anti)$ | $A$ | $-\varphi, \frac{1}{\varphi}, +1, -1$ | 1 | (2,2) |

Table 3.2

**Remark 3.2.** *Gram matrices are formed by defining conjugates such as*

$\bar{Z}^{\mathbb{CH}} = a - b_1h^1 + b_2h^2 - b_3h^3\epsilon\mathbb{H}$, (Complex- Hiperbolik- Conjugate) , $\bar{Z}^{\mathbb{H}A} = -a + b_1h^1 + b_2h^2 - b_3h^3\epsilon\mathbb{H}$ $(Hiperbolik - Anti -\ Conjugate)$, $\bar{Z}^{\mathbb{C}A} = -a - b_1h^1 + b_2h^2 + b_3h^3\epsilon\mathbb{H}$, $\bar{Z}^{D} = a + b_1h^1 - b_2h^2 + b_3h^3\epsilon\mathbb{H}$, $(Dual\ Conjugate)$, *... and obtaining similar eigenvalues.*

A Gram matrix can be constructed for each conjugate. The golden ratio can be reached in the eigenvalues. From the definition of G, the $h^1 - h^2$ block is found in the form $[\mp 1, \mp 1; \mp 1, 0]$. The characteristic equation of such a matrix is $\lambda^2 \pm \lambda - 1 = 0$, that is, the golden ratio equation.

### 3.2. Gram Matrix–Fibonacci Connection- Golden Ratio

**Definition 3.2. (Fibonacci Sequences)** The sequences defined by the initial conditions

$F_0 = 0,\ F_1 = 1$ and

$$F_n = F_{n-1} + F_{n-2}, \qquad n \geq 2$$

recursion is called the Fibonacci numbers [11, 12].

**Definition 3.3. (Lucas Sequences)** The sequence defined by the initial conditions

$L_0 = 2,\ L_1 = 1$ and

$$L_n = L_{n-1} + L_{n-2}, \qquad n > 1$$

recursion is called the Lucas sequences [11,12].

**Teorem 3.3. (Binet Formula)** For Golden ratio $\varphi = \frac{1+\sqrt{5}}{2}$ and its conjugate $-\frac{1}{\varphi} = \frac{1-\sqrt{5}}{2} = \psi$ are

$$F_n = \frac{\varphi^n - \psi^n}{\sqrt{5}}, \qquad L_n = \varphi^n + \psi^n.$$

The eigenvalues of the Gram matrices were $\varphi$ and $-\frac{1}{\varphi}$. Therefore,

$$F_n = \frac{\varphi^n - \psi^n}{\sqrt{5}} \qquad (3.15)$$

is written. Where $-\frac{1}{\varphi} = \frac{1-\sqrt{5}}{2} = \psi$. This establishes the link between Gram matrices and Fibonacci.

**Theorem 3.4. (Fibonacci- Gram Idendity)** For the $h^1 - h^2$subblok matrices $A$ of the gram matrices $G_O,\ G_{\mathbb{H}}$, and $G_A$:

$$A^n = \begin{bmatrix} F_{n+1} & -F_n \\ -F_n & F_{n-1} \end{bmatrix},\ n \geq 1\ integer, A = \begin{bmatrix} 1 & -1 \\ -1 & 0 \end{bmatrix}. \qquad (3.16)$$

**Proof.** Proof is done inductively. For $n = 1$,

$$A = \begin{bmatrix} 1 & -1 \\ -1 & 0 \end{bmatrix} = \begin{bmatrix} F_2 & -F_1 \\ -F_1 & F_0 \end{bmatrix},\ F_1 = 1,\ F_0 = 0,$$

$$A^2 = \begin{bmatrix} 1 & -1 \\ -1 & 0 \end{bmatrix} \cdot \begin{bmatrix} 1 & -1 \\ -1 & 0 \end{bmatrix} = \begin{bmatrix} 2 & -1 \\ -1 & 1 \end{bmatrix} = \begin{bmatrix} F_3 & -F_2 \\ -F_2 & F_1 \end{bmatrix}, F_2 = 1, F_3 = 2,$$

$$A^3 = \begin{bmatrix} 3 & -2 \\ -2 & 1 \end{bmatrix} = \begin{bmatrix} F_4 & -F_3 \\ -F_3 & F_2 \end{bmatrix}, F_3 = 2, F_4 = 3,$$

Let's assume that it is true for $n$, $A^{n+1} = A^n \cdot A$

And so on we write

$$A^n = F_{n+1} \begin{bmatrix} 1 & 0 \\ 0 & 0 \end{bmatrix} - F_n \begin{bmatrix} 0 & 1 \\ 1 & 0 \end{bmatrix} + F_{n-1} \begin{bmatrix} 0 & 0 \\ 0 & 1 \end{bmatrix}. \qquad (3.17)$$

**Corollary 3.5.**

- **i.** $Det(A^n) = F_{n-1} \cdot F_{n+1} - F_n^2 = (-1)^n$ , $n \geq 1$ (Cassini – Simson Identity)
- **ii.** $tr(A^n) = F_{n-1} + F_{n+1} = L_n$ , $n \geq 1$

**Theorem 3.5. (B subblok Fibonacci- Gram Idendity)** For the $h^1 - h^2$subblok matrices $B$ of the gram matrices $G_O$, $G_{\mathbb{H}}$, and $G_A$:

$$B^n = (-1)^n \begin{bmatrix} F_{n-1} & F_n \\ F_n & F_{n+1} \end{bmatrix},\ n \geq 1\ integer, B = \begin{bmatrix} -1 & 1 \\ 1 & 0 \end{bmatrix}. \qquad (3.18)$$

**Proof:** The relation $B = -JA\,J^{-1}$ ($J$ is the appropriate permutation matrix) and is proven by induction by showing that $A$ and $B$ are similar matrices.

**Corollary 3.6.**

- **i.** $Det(B^n) = (-1)^n$ , $n \geq 1$ (Cassini – Simson Identity)
- **ii.** $tr(B) = (-1)^n L_n$ , $n \geq 1$

**Corollary 3.7.**

For matrix $A^n = \begin{bmatrix} F_{n+1} & -F_n \\ -F_n & F_{n-1} \end{bmatrix}$ then, we write $A^n = \frac{L_n}{2} I + \frac{L_{n-1}+L_{n+1}}{10} \begin{bmatrix} 1 & -2 \\ -2 & -1 \end{bmatrix}$.

**Remark 3.3.** *These results in this section (Fibonacci–Gram identity, matrix proof of Cassini identity, Lucas–trace connection and B_block) arise naturally from the enlarged perspective provided by the 4×4 left-multiplication matrix representation.*

### 3.3. Lie Bracket, Hybrid Lie Algebra and Tensorial Structure Constants

Since $\mathbb{H}$ is an associative algebra, $(\mathbb{H},[\cdot,\cdot])$ is a Lie algebra with the commutator bracket $[X,Y]=XY-YX$. According to Table 1, commutator for $\{h^1,h^2,h^2\}$ created as follows :

| $[e_i,e_j]$ | $e_1=h^1$ | $e_2=h^2$ | $e_3=h^3$ |
|---|---|---|---|
| $e_1=h^1$ | 0 | $-e_3=-h^3$ | $e_1+e_2=h^1+h^2$ |
| $e_2=h^2$ | $e_3=h^3$ | 0 | $-e_2=-h^2$ |
| $e_3=h^3$ | $-e_1-e_2=-h^1-h^2$ | $e_2=h^2$ | 0 |

Table 3.3 (Lie Bracket Table, normalized by dividing by 2)

**Theorem 3.6.** $(\mathbb{H},[\cdot,\cdot])$ is a Lie algebra[2] with the commutator bracket $[h^i,h^j]=h^ih^j-h^jh^i$.

**Proof:**

- **i.** Bilinearity: $[\cdot,\cdot]$ is bilinear over $\mathbb{R}$,
- **ii.** Anti Symmetry: $[h^i,h^j]=-[h^j,h^i]$, it seen from table.
- **iii.** Jacobi Idendity: Since commutator is bilinear and ant,symmetric, it suffices to verify the Jacobi identity for basis elements $\{h^i, i=1,2,3\}$. Then we find $\sigma_{i,j,k}\left[h^i,[h^j,h^k]\right]=0$, where $\sigma_{i,j,k}$ is circular sums. Moreover, the identity element 1 is central; $[1,h^i]=0$ for all $h^i\epsilon\mathbb{H}$ ■

For any two basis elements of the Lie algebra $\{e_i=h^i, i=1,2,3\}$ the Lie commutator of those two basis elements can be written in the form

$$[e_i,e_j]=C_{ij}^k e_k \qquad (3.19)$$

Here, the coefficients $C_{ij}^k$ are called structural constants. These constants completely define the structure of the Lie algebra. Then we write

$$[h^i,h^j]=\sum_k C_{ij}^k h^k. \qquad (3.20)$$

Let $C_\alpha=(C_{\alpha\beta}^\gamma)$ be the matrix of structural constants. Therefore, according to the Table 2. We find $C_1=I$,

$$C_2=\begin{bmatrix}0&-1&1&0\\1&0&0&1\\0&0&0&1\\0&0&-1&0\end{bmatrix}, C_3=\begin{bmatrix}0&1&0&0\\0&0&0&0\\1&0&0&-1\\0&1&0&0\end{bmatrix}, C_4=\begin{bmatrix}0&0&0&1\\0&-1&0&0\\0&-1&1&0\\1&0&0&0\end{bmatrix} \qquad (3.21)$$

$[C_\alpha,C_\beta]=C_\alpha C_\beta-C_\beta C_\alpha$ commutators give matrix representations of $[h^i,h^j]$.

**Definition 3.4.** *Let $G$ is a Lie algebra of finite dimension on a field $F$. For $\forall X,Y\in G$,*

$$K(X,Y)=tr(ad_X{}^\circ ad_Y) \qquad (3.22)$$

*is called the Killing form on* $G$. *Where* $ad_X(Z) = [X, Z]$ *and* $ad_X$ *is adjoint matrix representation of* $X$ [6, 8 ,10].

According to the Table 2, we write

$$[h^1, h^2] = -h^3\ , [h^2, h^1] = h^3, [h^1, h^3] = h^1 + h^2, [h^3, h^1] = -h^1 - h^2$$
$$[h^2, h^3] = -h^2, [h^3, h^2] = h^2. \qquad (3.23)$$

Since $ad_{h^i}(h^j) = [h^i, h^j]$, we find that

$$[h^1, h^1] = 0, [h^1, h^2] = -h^3, [h^1, h^3] = h^1 + h^2$$

$$(ad_{h^1}) = \begin{bmatrix} 0 & 0 & 1 \\ 0 & 0 & 1 \\ 0 & -1 & 0 \end{bmatrix}, (ad_{h^2}) = \begin{bmatrix} 0 & 0 & 0 \\ 0 & 0 & -1 \\ 1 & 0 & 0 \end{bmatrix}, (ad_{h^3}) = \begin{bmatrix} -1 & 0 & 0 \\ -1 & 1 & 0 \\ 0 & 0 & 0 \end{bmatrix} \qquad (3.24)$$

According to the definition of Killing form, we find that

$$K(h^1, h^1) = tr(ad_{h^1} \circ ad_{h^1}) = -2, \quad K(h^2, h^2) = tr(ad_{h^2} \circ ad_{h^2}) = 2,$$
$$K(h^3, h^3) = tr(ad_{h^3} \circ ad_{h^3}) = 0, \qquad K(h^1, h^2) = tr(ad_{h^1} \circ ad_{h^2}) = 2,$$
$$K(h^1, h^3) = tr(ad_{h^1} \circ ad_{h^3}) = 0, \qquad K(h^2, h^3) = tr(ad_{h^2} \circ ad_{h^3}) = 0.$$

Thus we write matrix of Killig form as:

$$K = \begin{bmatrix} -2 & 2 & 0 \\ 2 & 0 & 0 \\ 0 & 0 & 2 \end{bmatrix}. \qquad (3.25)$$

Since $Det(K) = -8 \neq 0$, we find

$$K^{-1} = \begin{bmatrix} 0 & \frac{1}{2} & 0 \\ \frac{1}{2} & \frac{1}{2} & 0 \\ 0 & 0 & \frac{1}{2} \end{bmatrix}. \qquad (3.26)$$

To present the geometric characterization of the Lie algebra in a more refined form, the Killing form $K$ is normalized by a factor of ½ to obtain the matri $B = \frac{1}{2}K$. This normalization preserves the semi-simplicity and signature properties of the algebra while revealing deep mathematical connections whitin the eigenvalue spectrum. The characteristic equation of matrix $B$ is given by:

$$Det(B - \lambda I) = (1 - \lambda)(\lambda^2 + \lambda - 1) = 0. \qquad (3.27)$$

The eigenvalue are as follows:

$$\lambda_1 = 1,\ \lambda_2 = \frac{-1+\sqrt{5}}{2} = \frac{1}{\varphi}, \lambda_3 = \frac{-1-\sqrt{5}}{2} = -\varphi, \qquad (3.28)$$

Where $\varphi = \frac{1+\sqrt{5}}{2}$ represent the **Golden Ratio[4].** The fact that the eigenvalues can be expresseding term of Golden Ratio and its reciprocal $(\frac{1}{\varphi})$ indicates that the dynamic structure of the hybrid Lie algebra exists in a state of natural balance and fractal harmony. This enhances the potential for the system to be utalized in physical models, such as crystal structures or quasi-periodic systems.

**Theorem 3.5.** *Let* $(\mathbb{H}, [\cdot,\cdot])$ *be the Lie algebra spanned by the hybrid units* $\{h^i, i = 1,2,3\}$. *Fort he Killing form* $K(h^i, h^j) = tr(ad_{h^i} \circ ad_{h^j})$ *defined on this algebra, the following asserions hold:*

**i.** *The form* $K(h^i, h^j)$ *is non degenerate* $Det\ (K(h^i, h^j)) \neq 0$, *and the algebra is semi-simple*

**ii.** *The transformation* $\emptyset: h^i \to -h^i$ *leaves the Killing form invariant:*

$$K\left(\emptyset(h^i), \emptyset(h^j)\right) = K(h^i, h^j), h^i, h^j \in \mathbb{H},$$

**iii.**

$$K(ad_{h^i} h^j, h^k) + K(h^j, ad_{h^k} h^i) = 0.$$

**Proof.**

**i.** Since the matrix of $K$ Killing form

$$K = \begin{bmatrix} -2 & 2 & 0 \\ 2 & 0 & 0 \\ 0 & 0 & 2 \end{bmatrix} \Rightarrow Det(K) = -8 \neq 0,$$

the algebra is semi-simple.

**ii.** Killing form is a symmetric bilinear form, tahat is $K(\alpha X, Y) = \alpha K(X, Y)$ and $K(X, \alpha Y) = \alpha K(X, Y)$. Than we write

$$K\left(\emptyset(h^i), \emptyset(h^j)\right) = K(-h^i, -h^j) = K(h^i, h^j).$$

**iii.**

$$K(ad_{h^i} h^j, h^k) = K([h^i, h^i], h^k) = K(C_{ij}^l h^l, h^k) = C_{ij}^l K(h^l, h^k)$$

$$K(h^j, ad_{h^k} h^i) = K(h^j, [h^k, h^i]) = K(h^j, C_{ki}^l h^l)) = C_{ki}^l K(h^j, h^l) = C_{ji}^l K(h^k, h^l)$$

$$j \leftrightarrow k$$

$$K(ad_{h^i} h^j, h^k) + K(h^j, ad_{h^k} h^i) = C_{ij}^l K(h^l, h^k) + C_{ji}^l K(h^k, h^l)$$

$$= C_{ij}^l K_{lk} + C_{ji}^l K_{kl}.$$

Since $K_{lk} = K_{kl}$ and according to download process of indices,

$$C_{ijk} = C_{ij}^l K_{lk}, \qquad C_{jik} = C_{ji}^l K_{kl}$$

we find $C_{ijk} = -C_{jik}$, then $C_{ijk} + C_{jik} = 0$ ■

The idendity $C_{ijk} + C_{jik} = 0$ directly systems from the total antisymmetry of the structure constants when the indices are lowered using the Killing metrics. This property ensures that the adjoint action $ad_{h^i}$ is represented bye a skew symmetric matrix relative to Killing form, confirming that the Killing form is a bi-invariant metric on the hybrid Lie algebra.

**Definition 3.5. (Killing-Conjugate Pair)** For the killing form $K$ and the transformation $\emptyset: h^i \to -h^i$ the pair $(K, \emptyset)$ is called the Killing- Conjugate pair.

The Killing-Conjugate pair $(K, \emptyset)$ characterizes the hybrid algebra as a symmetric semi Riemannian manifold.

Since the real part $a \in \mathbb{R}$ lies in the center of the algebra $\mathbb{H}$ and does not contribute to the Lie bracket operations, we restrict our analysis to the pure hybrid subspace spanned by the basis $\{h^i, i = 1,2,3\}$.

**Definition 3.6. (Hybrid-Killing Metric)** Let $\mathbb{H}_0 = Span\{h^i, i = 1,2,3\}$ be the pure hybrid subspace of the Lie algebra $\mathbb{H}$. The unified metric

$$g_{ij} = Re(h^i \cdot \bar{h}^j) + \alpha K_{ij}, \alpha \in \mathbb{R} \quad (3.29)$$

is called Hybrid- Killing Metric on $\mathbb{H}_0$.

**Theorem 3.6.** Let $\mathbb{H}_0 = Span\{h^i, i = 1,2,3\}$ be the pure hybrid subspace of the Lie algebra $\mathbb{H}$. The unified metric $g_{ij} = Re(h^i \cdot \bar{h}^j) + \alpha K_{ij}$ is invariant under the conjugation $\emptyset: h^i \to -h^i$.

**Proof.** Let, $Re(h^i \cdot \bar{h}^j) = A_{ij}$ . Then, we write $g_{ij} = A_{ij} + \alpha K_{ij}$. we shown that

$$g_{ij}\left(\emptyset(h^i), \emptyset(h^j)\right) = g_{ij}(h^i, h^j).$$

By the bilinearity of the unified metric, the action of $g_{ij}$ on the conjugated basis elements $\emptyset(h^i)$ can be decomposed into the sum of the conjugated hybrid product and the conjugated Killing form. Then, we write

$$\begin{aligned} g_{ij}\left(\emptyset(h^i), \emptyset(h^j)\right) &= (A_{ij} + \alpha K_{ij})\left(\emptyset(h^i), \emptyset(h^j)\right) \\ &= A_{ij}\left(\emptyset(h^i), \emptyset(h^j)\right) + \alpha K_{ij}\left(\emptyset(h^i), \emptyset(h^j)\right). \end{aligned}$$

Since $\emptyset(h^i) = -h^i$, $\emptyset(h^j) = -h^j$ and $\overline{\emptyset(h^j)} = \overline{-h^j} = -\bar{h}^j$, then we write

$$Re(\emptyset(h^i) \cdot \overline{\emptyset(h^j)}) = Re(-h^i \cdot -\bar{h}^j) = Re(h^i \cdot \bar{h}^j) = A_{ij}.$$

Since

$$[\emptyset(h^i), \emptyset(h^j)] = [-h^i, -h^j] = [h^i, h^j], \qquad ad_{\emptyset(h^i)} = ad_{-h^i} = -ad_{h^i}$$

We find

$$K\left(\emptyset(h^i), \emptyset(h^j)\right) = tr\left(-ad_{h^i} \circ -ad_{h^j}\right) = tr\left(ad_{h^i} \circ ad_{h^j}\right) = K(h^i, h^j) = K_{ij} \blacksquare$$

Since, $A_{11} = 1$, $A_{22} = 0$, $A_{33} = -1$, $K_{11} = -2$, $K_{22} = 0$, $K_{33} = 2$ we find

$$g = \begin{bmatrix} 1 - 2\alpha & 0 & 0 \\ 0 & 0 & 0 \\ 0 & 0 & 2\alpha - 1 \end{bmatrix}. \quad (3.30)$$

The Hybrid-Killing metric $g_{ij}$ exhibits a unique topolojical phase transition governed by the coupling parameter $\alpha$. By analazing the metric signature, we identify $\alpha = 1/2$ as the critical threshold where the algebraic curvature of the Killing form perfectly counterbalance the hybrid inner product.

**Corollary 3.5.** For $\alpha < 1/2$, the metric maintains a Lorenzian-like character, where the causal structure is dominated by the hybrid real part. For $\alpha = 1/2$, at this singularity, the manifold becomes entirely isotropic (null). The cone of light geometrically "collapses" into a point or a plane. And, for $\alpha > 1/2$, the signature flips, signifies that the Lie algebra structure takes precedence. The orientation of space is completely reversed. Time-like orientations and space-like orientations switch places.

**Remark 3.4.** Istoropic cones (null cones) constitute the set of vectors for which a metric vanishes $g(X,X) = 0$. In hybrid space, these cones physically represent "critical directions that move at the speed of light or maintain a specific energy balance"

With respect to metric $g_{ij} = Re(h^i \cdot \bar{h}^j) + \alpha K_{ij}$ , $\alpha \in \mathbb{R}$, since $g_{11} = 1 - 2\alpha, g_{22} = 0, g_{33} = 2\alpha - 1$ we write $g_{11}(b_1)^2 + g_{22}(b_2)^2 + g_{33}(b_3)^2 = 0$. Then we find

$$(1-2\alpha)(b_1)^2 + (2\alpha - 1)(b_3)^2 = 0$$
$$(1-2\alpha)((b_1)^2 - (b_3)^2) = 0.$$

If $\alpha \neq 1/2$: The equation of the cone becomes $(b_1)^2 - (b_3)^2 = 0$, i.e. $b_1 = \mp b_3$. In hybrid space, this is a "slab or a set of V-shaped planes extending along the $h^2$ axis."
If $\alpha = 1/2$ : The metric becomes completely degenerate (all coefficients vanish). In this case, the entire space becomes an 'isotropic universe'; the concept of distance disappears."

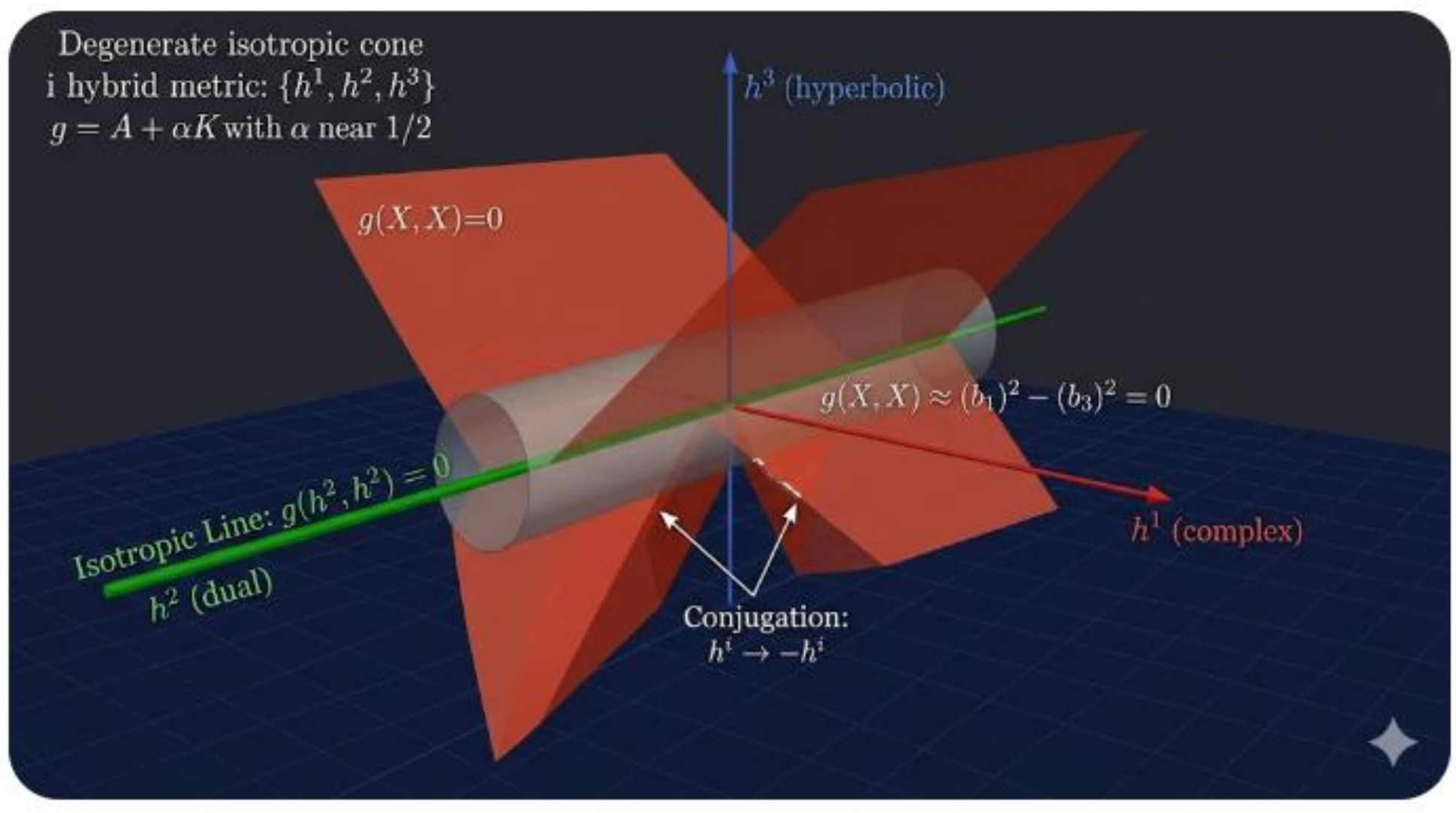

Figure 1 (was generated by AI)

## 4. Hybrid - Type Decomposition

Multiplying two compleks number $a + bI$ and $c + dI$ is straightforward

$$(a,b)(c,d) = (ac - bd, ad + bc).$$

For two quaternions , $bI$ and $dI$ become the 3-vectors $B$ abd $D$, where $B = xI + yJ + zK$ and similarly $D$. Multiplication of quaternions is like complex numbers, but with the addition of the cross product [ 7]:

$$(a,\vec{B})(c,\vec{D}) = (ac - \vec{B}\cdot\vec{D}, a\vec{D} + c\vec{B} + \vec{B}\times\vec{D}).$$

If a is the operator d/dt, and B is the del operator, or d/dx I + d/dy J + d/dz K (all partial derivatives), then these operators act on the scalar function c and the 3−vector function D in the following manner:

$$\left(\frac{d}{dt},\vec{\nabla}\right)\left(c,\vec{D}\right)=\left(\frac{dc}{dt}-\vec{\nabla}\cdot\vec{D},\frac{d\vec{D}}{dt}+\vec{\nabla}c+\vec{\nabla}\times\vec{D}\right).$$

This one quaternion contains the time derivatives of the scalar and 3−vector functions, along with the divergence, the gradient and the curl[7].

$\vec{V_1}$ and $\vec{V_2}$, hybrid numbers are like complex numbers,

$$Z_1\cdot Z_2=\left(a_1a_2-S\left(\vec{V_1},\vec{V_2}\right),\quad a_1\cdot\vec{V_2}+a_2\cdot\vec{V_1}+\vec{V_1}\times_{\mathbb{H}}\vec{V_2}\right).$$

We will call this hybrid-type decomposition.

Let two hybrid numbers, $Z_1,Z_2\epsilon\mathbb{H}$, $Z_1=a_1+b_1h^1+c_1h^2+d_1h^3$, $Z_2=a_2+b_2h^1+c_2h^2+d_2h^3$, $a_i,c_i,b_i,d_j\epsilon\mathbb{R},i,j=1,2$. Let us, $Z_1=\left(a_1,\vec{V_1}\right)$, $Z_2=\left(a_2,\vec{V_2}\right)$.

Here $\vec{V_1}=b_1h^1+c_1h^2+d_1h^3$ and $\vec{V_2}=b_2h^1+c_2h^2+d_2h^3$.

$$Z_1\cdot Z_2=\underbrace{(a_1a_2-b_1b_2+b_1c_2+b_2c_1+d_1d_2)}_{Real\ part}+\underbrace{(a_1b_2+a_2b_1+b_1d_2-b_2d_1)}_{h^1\ component}h^1$$
$$+\underbrace{(a_1c_2+a_2c_1+b_1d_2-b_2d_1-c_1d_2+c_2d_1)h^2}_{h^2\ component}+\underbrace{(a_1d_2+a_2d_1-b_1c_2+b_2c_1)}_{h^3\ component}h^3.\quad(4.1)$$

If we seperate two part (4.1) as

$$a_1\cdot\vec{V_2}+a_2\cdot\vec{V_1}=(a_1b_2+a_2b_1)h^1+(a_1c_2+a_2c_1)h^2+(a_1d_2+a_2d_1)h^3,$$
$$(b_1d_2-b_2d_1)h^1+(b_1d_2-b_2d_1-c_1d_2+c_2d_1)h^2+(-b_1c_2+b_2c_1)h^3.$$

Now, we write real part of (4.1)

$$Re(Z_1\cdot Z_2)=a_1a_2-S\left(\vec{V_1},\vec{V_2}\right).$$

**Definition 4.1. (Hybrid scalar product**) The $S\left(\vec{V_1},\vec{V_2}\right)$ is called hybrid scalar product ,

$$S\left(\vec{V_1},\vec{V_2}\right)=b_1b_2-b_1c_2-c_1b_2-d_1d_2.\qquad(4.2)$$

**Definition 4.2.** The space of vector parts of hybrid numbers is defined as

$$\vec{\mathbb{H}}_0=\left\{\vec{V}=bh^1+ch^2+dh^3:b,c,d\in\mathbb{R}\right\}.$$

The hybrid vector product on this space is defined as follows.

**Definition 4.3. (Hybrid vectorial product)**

The $(b_1d_2-b_2d_1)h^1+(b_1d_2-b_2d_1-c_1d_2+c_2d_1)h^2+(-b_1c_2+b_2c_1)h^3$ term of equation (4.1) is called the vector product of $\vec{V_1}$ and $\vec{V_2}$ and is denoted by $\vec{V_1}\times_{\mathbb{H}}\vec{V_2}$.

$$\overrightarrow{V_1} \times_{\overrightarrow{\mathbb{H}}} \overrightarrow{V_2} = (b_1 d_2 - b_2 d_1)h^1 + (b_1 d_2 - b_2 d_1 - c_1 d_2 + c_2 d_1)h^2 + (-b_1 c_2 + b_2 c_1)h^3. \quad (4.3)$$

Then, we write,

$$Z_1 \cdot Z_2 = \left(a_1 a_2 - S(\overrightarrow{V_1}, \overrightarrow{V_2}),\ a_1 \cdot \overrightarrow{V_2} + a_2 \cdot \overrightarrow{V_1} + \overrightarrow{V_1} \times_{\overrightarrow{\mathbb{H}}} \overrightarrow{V_2}\right). \quad (4.4)$$

**Corolary 4.1**

- **i.** $S(\overrightarrow{V_1}, \overrightarrow{V_2}) = S(\overrightarrow{V_2}, \overrightarrow{V_1})$,
- **ii.** $S(\overrightarrow{V_1}, \overrightarrow{V_1}) = b^2 - 2bc - d^2 = (b-c)^2 - c^2 - d^2,\ \vec{V} = bh^1 + ch^2 + dh^3.$
- **iii.** $Re\,(Z_1 \cdot Z_2) = a_1 a_2 - S(\overrightarrow{V_1}, \overrightarrow{V_2})$.
- **iv.** The hybrid scalar product is pseudo metrics.

From (4.3), we write $S\left(\vec{V}, \vec{V}\right) = (b-c)^2 - c^2 - d^2$ and has signature (1,2). Therefore $S$ is not an inner product but a **pseudo- metric bilineer form.**

$$\overrightarrow{V_1} \times_{\mathbb{H}} \overrightarrow{V_2} = \begin{vmatrix} b_1 & b_2 \\ d_1 & d_2 \end{vmatrix} h^1 + \left( \begin{vmatrix} b_1 & b_2 \\ d_1 & d_2 \end{vmatrix} - \begin{vmatrix} c_1 & c_2 \\ d_1 & d_2 \end{vmatrix} \right) h^2 + \begin{vmatrix} b_2 & b_1 \\ c_2 & c_1 \end{vmatrix} h^3. \quad (4.5)$$

**Corolary 4.2.**

- **i.** For any hybrid number $Z = a + \vec{V} = a + bh^1 + ch^2 + dh^3$, we have $\vec{V} \times_{\overrightarrow{\mathbb{H}}} \vec{V} = \vec{O}$.
- **ii.** $\overrightarrow{V_1} \times_{\overrightarrow{\mathbb{H}}} \overrightarrow{V_2} = -(\overrightarrow{V_2} \times_{\overrightarrow{\mathbb{H}}} \overrightarrow{V_1})$.
- **iii.** Jacoby idendity:
  $$\overrightarrow{V_1} \times_{\overrightarrow{\mathbb{H}}} (\overrightarrow{V_2} \times_{\overrightarrow{\mathbb{H}}} \overrightarrow{V_3}) + \overrightarrow{V_2} \times_{\overrightarrow{\mathbb{H}}} (\overrightarrow{V_3} \times_{\overrightarrow{\mathbb{H}}} \overrightarrow{V_1}) + \overrightarrow{V_3} \times_{\overrightarrow{\mathbb{H}}} (\overrightarrow{V_1} \times_{\overrightarrow{\mathbb{H}}} \overrightarrow{V_2}) = \vec{O} \quad (4.6)$$

**Theorem 4.1. (Hybrid Lie Algebra)** $(\overrightarrow{\mathbb{H}},\ \times_{\overrightarrow{\mathbb{H}}})$ is a Lie algebra.

**Proof.** We verify the tree axioms of a Lie algebra

- **i.** Bilinearity: Each componentt of the hybrid vektör product is a linear combination of 2 x 2 determinants in coordinates. Since determinants are bilinear, $\times_{\overrightarrow{\mathbb{H}}}$ is linear in eachargument.
- **ii.** Anti symmetry: For all $\overrightarrow{V_1},\ \overrightarrow{V_2} \in \overrightarrow{\mathbb{H}},\ \overrightarrow{V_1} \times_{\overrightarrow{\mathbb{H}}} \overrightarrow{V_2} = -(\overrightarrow{V_2} \times_{\overrightarrow{\mathbb{H}}} \overrightarrow{V_1})$.
- **iii.** Jacobi Idendity: For $\overrightarrow{V_1},\ \overrightarrow{V_2}, \overrightarrow{V_3} \in \overrightarrow{\mathbb{H}}$,
  $$\overrightarrow{V_1} \times_{\overrightarrow{\mathbb{H}}} (\overrightarrow{V_2} \times_{\overrightarrow{\mathbb{H}}} \overrightarrow{V_3}) + \overrightarrow{V_2} \times_{\overrightarrow{\mathbb{H}}} (\overrightarrow{V_3} \times_{\overrightarrow{\mathbb{H}}} \overrightarrow{V_1}) + \overrightarrow{V_3} \times_{\overrightarrow{\mathbb{H}}} (\overrightarrow{V_1} \times_{\overrightarrow{\mathbb{H}}} \overrightarrow{V_2}) = \vec{O}.$$

We write $\overrightarrow{V_2} \times_{\overrightarrow{\mathbb{H}}} \overrightarrow{V_3} = \alpha h^1 + \beta h^2 + \gamma h^3$, where

$$\alpha = (b_2 d_3 - b_3 d_2), \qquad \beta = \alpha - (c_2 d_3 - c_3 d_2), \qquad \gamma = (b_3 c_2 - b_2 c_3)$$

The components of $\overrightarrow{V_1} \times_{\overrightarrow{\mathbb{H}}} (\overrightarrow{V_2} \times_{\overrightarrow{\mathbb{H}}} \overrightarrow{V_3})$ are summed with the remaining two terms obtained by the cyclic permutation $1 \to 2 \to 3 \to 1$.

**Remark 4.1.** The space $\overrightarrow{\mathbb{H}}$ equipped with hybrid vector product $\times_{\overrightarrow{\mathbb{H}}}$ forms a Lie algebra that is analogous to the standart cross product on $\mathbb{R}^3$, yet reflects the distinctive algebraic structure of hybrid numbers. This result shows that the vector part of hybrid numbers is not marely a vector space, but carries a rich algebraic geometric intrinsive to the hybrid number system.

From (4.5) we write $Comp\ h^2 = Comph^1 - (c_1d_2 - c_2d_1)$. Than we rewrite hiybrid vectorial product as:

$$\overrightarrow{V_1} \times_{\mathbb{H}} \overrightarrow{V_2} = \underbrace{(b_1d_2 - b_2d_1)}_{\Delta_{bd}} h^1 + (\Delta_{bd} - \Delta_{cd})h^2 + \underbrace{(b_2c_1 - b_1c_2)}_{\Delta_{bc}} h^3$$
$$= \Delta_{bd}h^1 + (\Delta_{bd} - \Delta_{cd})h^2 + \Delta_{bc}h^3. \qquad (4.7)$$

Where,

$$\Delta_{bd}= \begin{vmatrix} b_1 & b_2 \\ d_1 & d_2 \end{vmatrix}, \Delta_{cd}= \begin{vmatrix} c_1 & c_2 \\ d_1 & d_2 \end{vmatrix}, \Delta_{bc}= \begin{vmatrix} b_2 & b_1 \\ c_2 & c_1 \end{vmatrix}. \qquad (4.8)$$

**Theorem 4.2 (Decomposition Theorem for Hybrid Numbers)**

$Z_1 = a_1 + \overrightarrow{V_1}$ and $Z_2 = a_2 + \overrightarrow{V_2}$ be two hybrid numbers with vector parts $\overrightarrow{V_k} = b_kh^1 + c_kh^2 + d_kh^3$. Then the multiplication of hybrid numbers admits the following decomposition for hybrid numbers:

$$Z_1 \cdot Z_2 = \left(a_1a_2 - S\left(\overrightarrow{V_1}, \overrightarrow{V_2}\right), \qquad a_1 \cdot \overrightarrow{V_2} + a_2 \cdot \overrightarrow{V_1} + \overrightarrow{V_1} \times_{\overrightarrow{\mathbb{H}}} \overrightarrow{V_2}\right).$$

In standars 3D vector algebra over real numbers, the cross product of two vector $\overrightarrow{V_1} \times \overrightarrow{V_2}$ produce a perpendicular vector. For hybrid numbers, we define a similar operation if we treat the imaginary components as a vector.

***Remark.*** *If $a$ is the operator $\frac{d}{dt}$, and $\overrightarrow{V_1}$ is the del operator, or $\frac{\partial}{\partial x}h^1 + \frac{\partial}{\partial y}h^2 + \frac{\partial}{\partial z}h^3$ (all partial derivatives), then these operators act on the function $c$ and the $3-$vector function $\overrightarrow{V_2}$ in the following manner:*

$$\left(\frac{d}{dt}, \vec{\nabla}\right)\left(c, \overrightarrow{V_2}\right) = \left(\frac{dc}{dt} - \vec{\nabla} \cdot \overrightarrow{V_2}, \frac{d\overrightarrow{V_2}}{dt} + \vec{\nabla}c + \vec{\nabla} \times \overrightarrow{V_2}\right)$$

*This contains the time derivatives of the scalar and 3−vector functions, along with the divergence, the gradient and the curl. This terminology leads us to physical applications.*